\documentclass{baustms}

\citesort

\theoremstyle{cupthm}
\newtheorem{thm}{Theorem}[section]

\newtheorem{lemma}[thm]{Lemma}
\theoremstyle{cupdefn}

\theoremstyle{cuprem}
\newtheorem{rem}[thm]{Remark}
\numberwithin{equation}{section}

\usepackage{verbatim}

\newcommand{\bR}{{\mathbb R}}
\newcommand{\bZ}{{\mathbb Z}}

\newcommand{\bfi}{{\mathbf i}}
\newcommand{\bfj}{{\mathbf j}}

\newcommand{\bfk}{{\mathbf k}}

\newcommand{\bfw}{{\mathbf w}}

\newcommand{\bfx}{{\mathbf x}}
\newcommand{\bfu}{{\mathbf u}}
\newcommand{\bfv}{{\mathbf v}}
\newcommand{\bfp}{{\mathbf p}}

\newcommand{\ctbfw}{\sigma({\mathbf w})}

\newcommand{\upM}{{\overline{M} }}
\newcommand{\lowM}{{\underline{M}}}
\newcommand{\ups}{{\overline{s} }}
\newcommand{\lows}{{\underline{s}}}

\newcommand{\tdM}{{\widetilde{M} }}
\newcommand{\tdu}{{\tilde{u}}}

\newcommand{\red}{\textcolor{red}}

\newcommand{\green}{\textcolor{green}}

\newcommand{\Var}{{\mathrm{Var}}}

\def\gep{\varepsilon}
\def\N{\mathcal{N}}

\begin{document}
\runningtitle{Box dimension of fractal interpolation surfaces}
\title{Box dimension of fractal interpolation surfaces with vertical scaling function}
\author[1]{Lai Jiang}
\address[1]{School of Mathematical Sciences, Zhejiang University, Hangzhou 310058, China\email{jianglai@zju.edu.cn}}

\authorheadline{Lai Jiang}


\support{The research was supported in part by the NSFC grants 12371089 and by ZJNSFC grant LY22A010023.}

\begin{abstract}
In this paper, we first present a simple lemma which allows us to estimate the box dimension of graphs of given functions by the associated oscillation sums and oscillation vectors.
Then we define vertical scaling matrices of generalized affine fractal interpolation surfaces (FISs). 
By using these matrices, we establish relationships between oscillation vectors of different levels, which enables us to obtain the box dimension of generalized affine FISs under certain constraints.

\end{abstract}

\classification{primary 28A80; secondary 41A30}
\keywords{fractal interpolation surfaces, box dimension, vertical scaling matrices}

\maketitle

\section{Introduction}
Fractal functions is one of the fundamental objects in the fractal geometry.
There are many works on fractal functions, including Weierstrass functions \cite{S18} and Tagaki functions \cite{AK11}.
In this paper, we will focus on a special class of fractal functions called fractal interpolation functions (FIFs), which was introduced by Barnsley \cite{B86}.
The graphs of FIFs are invariant sets of certain iterated function systems (IFSs).
We remark that both Weierstrass functions and Tagaki functions are FIFs.
Moreover, in the case that every map in the corresponding IFS is affine, the graph of the FIF is an self-affine set, and the function is usually called an affine FIF. 
There are many works on affine FIFs, such as \cite{BRS20,B89,HM86}.

FIFs are very useful in the field of interpolation and fitting.
In order to fit bivariate functions, it is natural to construct fractal interpolation surfaces (FISs).
In \cite{M90}, Massopust introduced FISs on triangles, where the interpolation points on the boundary are required to be coplanar.
Metzler and Yun \cite{MY10} presented a method to generate FISs on an arbitrary data set over rectangular grids. 
This method was developed by Ruan and Xu \cite{RX15}, where they introduced a general framework to construct FISs on rectangular grids. 
They also defined bilinear FISs in that paper.
There are other methods to construct FISs, including \cite{BD07,F08,VV20,YCO15}.


Finding dimensions of fractal sets, including box and Hausdorff dimensions, is one of the most concerned problems in fractal geometry.
Hardin and Massopust \cite{HM86} and Barnsley \emph{et al}. \cite{BEHM89} obtained the box dimension of the graphs of affine FIFs. 
There are also other works on the box dimension of the graphs of FIFs, please see \cite{BM15,RSY09,RXY21}.
Recently, Jiang and Ruan \cite{JR23,JR24} studied the box dimension of generalized affine FIFs, where the vertical scaling functions satisfy Lipschitz conditions.

There has been some works on estimating the box dimension of FISs, see \cite{BD07,F08,VV20,YCO15}.
In \cite{F08}, Feng obtained the exact box dimension of FISs with vertical scaling factors.
It is very challenging to obtain the exact box dimension of FISs with vertical scaling functions.
So far as we know, there are few results on this direction. 
In the work of Kong \emph{et al.} \cite{KRZ18}, they required that the vertical scaling functions are bilinear functions with a uniform sum.
In this paper, we will study this problem under the condition that vertical scaling functions are Lipschitz.
We mainly follow the spirit of \cite{JR23,JR24}. 
We first define vertical scaling matrices and oscillation vectors.
Then we use vertical scaling matrices to establish relationships between oscillation vectors of different levels, which provides a way to obtain the exact box dimension under certain constraints.
We remark that in order to prove our main results (see the proof of Lemma \ref{lem:infty} and Theorem \ref{th:lower}),
we need to normalize oscillation vectors and vertical matrices, which is different to one-dimensional case.

The paper is organized as follows. 
In Section 2, we recall some basic notations and definitions.
In Section 3, we define vertical scaling matrices and oscillation vectors.
In Section 4, we prove main results and obtain the box dimension of FISs under certain constraints.

\section{Preliminaries}

\subsection{Definition of fractal interpolation surfaces.}

We follow the basic idea of \cite{MY10} and \cite{RX15}.
Given a data set $\big\{(x_{i}, y_{j},z_{ij})\in \mathbb{R}^{3}:\, i\in\{0,1,\ldots,N\},j\in\{0,1,\ldots,M\} \big\}$ with
$$x_{0}<x_{1}<\cdots<x_{N}, \quad y_{0}<y_{1}<\cdots<y_{M},$$
where $N,M\geq 2$ are positive integers. 

Let $I:=[x_0,x_N]$ and $J:=[y_0,y_M]$. Write $D:=I \times J$.
For any $1 \leq i \leq {N}$, we write $I_{i}:=[x_{i-1},x_{i}]$. Let $u_i:\, I\to I_i$ be the affine function satisfying
\begin{align}
	&  u_i(x_0)=x_{i-1}, \quad u_i(x_N)=x_i, 	\quad\quad \mbox{if $i$ is odd, and} \label{eq:ui-def-1} \\
	&  u_i(x_0)=x_{i},  \quad\quad u_i(x_N)=x_{i-1},	\quad \mbox{if $i$ is even.} \label{eq:ui-def-2}
\end{align}
Clearly, this implies that $u_1(x_0)=x_0$, $u_1(x_N)=u_2(x_N)=x_1$, $u_2(x_0)=u_3(x_0)=x_2$ and so on. 

Similarly, for any $1 \leq j \leq M$, we write $J_{j}:=[y_{j-1},y_{j}]$. Let $v_j:\, J\to J_j$ be the affine function satisfying
\begin{align}
&  v_j(y_0)=y_{j-1}, \quad v_j(y_M)=y_j, \quad\quad \mbox{if $j$ is odd, and} \label{eq:vj-def-1}\\
&  v_j(y_0)=y_{j}, \quad\quad v_j(y_M)=y_{j-1}, \quad \mbox{if $j$ is even.} \label{eq:vj-def-2}
\end{align}

We regard the digit set $\Sigma=\{1,2,\cdots,N\} \times  \{1,2,\cdots,M\}$ as the index set.
For any $w=(i,j) \in \Sigma$,
write $D_w:=I_i\times J_j$ and let $L_w:D \to D_w$ be a continuous mapping defined by
\[
	L_{w}(x,y)=\big(u_i(x),v_j(y)\big).
\]

Let $g:\, D\to \bR$ be a continuous function satisfying
\begin{equation*}
 g(x_i,y_j)=z_{ij}, \quad \forall \, (i,j)\in \{0,N\}\times \{0,M\}.
\end{equation*}
Let $h:\, D\to \bR$ be a continuous function satisfying
\begin{equation*}
 h(x_i,y_j)=z_{ij}, \quad \forall \, (i,j)\in \{0,1,\ldots,N\} \times \{0,1,\ldots,M\}.
\end{equation*}
Let $S:\, D\to \bR$ be a continuous function satisfying $|S(x,y)|<1$ for all $(x,y)\in D$.


For any $w\in \Sigma$, we define $F_{w}:\, D\times \bR\to \bR$ by
\begin{equation}\label{eq:Fij-def-MetYun}
  F_{w}(x,y,z)=S\big(L_w(x,y)\big)\big(z-g(x,y)\big)+h\big(L_w(x,y)\big).
\end{equation}

Now, for any $w=(i,j) \in \Sigma$, we define $\Psi_w:\, D\times \bR\to D_w \times \bR$ by
\begin{equation}\label{eq:Wij-def}
  \Psi_{w}(x,y,z)=\big(u_i(x),v_j(y), F_w(x,y,z)\big).
\end{equation}
Then $\{D\times \bR, \Psi_w: w \in \Sigma \}$ is an \emph{iterated function system} (IFS).

\begin{thm}[\cite{MY10,RX15}]
Let $\{ D \times \bR, \Psi_{w}:  w \in \Sigma \}$ be the IFS defined in \eqref{eq:Wij-def}.
Then there exists a unique continuous function $f:\,D\to \bR$ such that $f(x_i,y_j)=z_{ij}$ for all $0 \leq i \leq N$ and $0 \leq j \leq M$ and $\Gamma f=\bigcup_{w \in \Sigma} \Psi_w (\Gamma f)$, where $\Gamma f=\big\{(x,y,f(x,y)):\, (x,y)\in D \big\}$ is the graph of $f$.  
We call $\Gamma f$ the \emph{generalized affine fractal interpolation surface} (FIS) and $f$ the \emph{generalized affine fractal interpolation function} (FIF) with respect to the IFS $\{ D \times \bR, \Psi_w: w \in \Sigma \}$.
\end{thm}

The function $S$ is called the \emph{vertical scaling function} of $\Gamma f$.
In the case that $g$ is a bilinear function, and both $h$ and $S$ are piecewise bilinear functions, Ruan and Xu \cite{RX15} called $\Gamma f$ a \emph{bilinear FIS}.

\subsection{Basic settings in this paper.}

Since the box dimension is one of basic dimensions in fractal geometry, it is quite natural to study the box dimension of FISs. In fact, Kong, Ruan and Zhang \cite{KRZ18} obtained the box dimension of bilinear FIS under the assumption that $N=M$, $x_i=y_i=i/N$ for all $0\leq i\leq N$, and the vertical scaling factors are steady and have uniform sum. We remark that this result was generalized by Liang and Ruan in \cite{LR21}, where they obtained the box dimension of bilinear recurrent FISs under certain constraints.

In the present paper, we will study the box dimension of generalized affine FISs. 

First,  we assume that $N=M \geq2$ and $|I|=|J|$. 
Moreover, we assume that $x_i-x_{i-1}=(x_N-x_0)/N$ and $y_i-y_{i-1}=(y_N-y_0)/N$ for each $1 \leq i \leq N$.
Under these assumptions, for any $w=(i,j)\in \Sigma$,
\begin{equation*}
	L_w(x,y)=\big( a_i x +b_i,c_j y+d_j\big) ,
\end{equation*}
where $a_i,b_i,c_j,d_j \in \mathbb{R}$ with $a_i=(-1)^{i+1}/N$, $c_j=(-1)^{j+1}/N$.

For convenience, we define a function on $D$ by
\begin{equation}\label{eq:q-def}
q(x,y)=h(x,y)-S(x,y)g\big(L_w^{-1}(x,y)\big), \quad (x,y) \in D_w.
\end{equation}
It is easy to check that $q$ is well-defined on $\bigcup_{w \in \Sigma}\partial D_w$ and $q$ is continuous on $D$, where $\partial D_w$ is the boundary of $D_w$. 
From \eqref{eq:q-def}, we can
 rewrite \eqref{eq:Fij-def-MetYun} as follows.
\begin{equation*}
	F_{w}(x,y,z)=S\big(L_w(x,y)\big) z+q\big(L_w(x,y)\big).
\end{equation*}
Furthermore, we assume that the following two conditions are satisfied.
\begin{itemize}
	\item[(A1)] $S$ is Lipschitz on $D$, i.e., there exists $\lambda_S>0$, such that $\big|S(t_1)-S(t_2)\big|<\lambda_S |t_1-t_2|$ for all $t_1,t_2 \in D$,
	\item[(A2)] $q$ is Lipschitz on $D$, i.e., there exists $\lambda_q>0$, such that $\big|q(t_1)-q(t_2)\big|<\lambda_q |t_1-t_2|$ for all $t_1,t_2 \in D$.
\end{itemize}

We remark that the function $q$ defined by \eqref{eq:q-def} is Lipschitz on $D$ if functions $S$, $g$ and $h$ are Lipschitz on $D$.

Under the above settings, the following notations are widely used.

\begin{itemize}
	\item[(1)] For $n \in \mathbb{Z}^+$, write $\Sigma^n := \{ \bfw =w_1 \cdots w_n : w_1,\ldots, w_n \in \Sigma\}$.
			 Let $\Sigma^0:=\{ \vartheta\}$, where $\vartheta$ denotes the \emph{empty word}. 
			 For $n \geq 0$ and $\bfw \in \Sigma^n$, we call $\bfw$ a \emph{word} with length $|\bfw|=n$.
	\item[(2)] Let $\Sigma^*:=\bigcup_{n=0}^\infty \Sigma^n$ be the collection of words with finite length.
	\item[(3)] For $n \geq 1$ and $\bfw=w_1  w_2 \cdots w_n \in \Sigma^n$, we define $\ctbfw:=w_2 \cdots w_n \in \Sigma^{n-1}$.
	\item[(4)] For any $n,m \geq 1$, $\bfp=p_1 \cdots p_n \in \Sigma^n$ and $\bfw=w_1 \cdots w_m \in \Sigma^m$, we write 
	\[\bfp\bfw:=p_1 \cdots p_n w_1 \cdots w_m \in \Sigma^{n+m}.\] 
	Write $\bfp \vartheta=\bfp=\vartheta \bfp$ for all $\bfp \in \Sigma^*$.
	\item[(5)] For $n \geq 1$ and $\bfw=w_1 \cdots w_n \in \Sigma^n$, let $L_\bfw:=L_{w_1}\circ \cdots \circ L_{w_n}$ and $D_{\bfw}= L_\bfw(D)$. Let $D_\vartheta=D$.
\end{itemize}

\subsection{Oscillation sums and box dimension.} 
Let $\varphi$ be a continuous function on $D$.
Given $n \in \mathbb{Z}^+$, 
we define \emph{oscillation sum} of $\varphi$ on $D$ with level $n$ by
\begin{align*}
  &O_n(\varphi,D)=  \sum_{\bfw \in \Sigma^{n} }O(\varphi,D_{\bfw}),
\end{align*}
where we use $O(\varphi,E)$ to denote the \emph{oscillation} of $\varphi$ on $E\subset D$, i.e.,
\[
	O(\varphi,E) = \sup \big\{\varphi(\bfx^\prime)-\varphi(\bfx^{\prime\prime}):\; \bfx^\prime,\bfx^{\prime\prime}\in E\big\}.
\]
For any $\bfp \in \Sigma^*$, we
define $O_n(\varphi,D_\bfp)= \sum_{\bfw \in \Sigma^{n} }O(\varphi,D_{\bfp\bfw})$.

Clearly, we have the following simple fact.
\begin{lemma}\label{lem:q-lips}
Let $\varphi$ be a Lipschitz continuous function on $D$, i.e. there exists $ \lambda \in \mathbb{R}$, $|\varphi(t_1)-\varphi(t_2)|< \lambda |t_1-t_2|$ for all $t_1,t_2 \in D$. 
Then for any $\bfp \in \Sigma^*$ and $n \in \mathbb{Z}^+$,
\begin{equation*}
	O_n(\varphi, D_\bfp)\leq  \sqrt{2}\lambda |I| N^{n-|\bfp|}.
\end{equation*}
\end{lemma}
\begin{proof}
Fix $n \in \mathbb{Z}^+$ and $\bfp \in \Sigma^*$. Then $O(\varphi,D_{\bfp \bfw}) \leq \lambda |D_{\bfp \bfw}|=\lambda \sqrt{2}|I|N^{-n-|\bfp|}$ for all $\bfw \in \Sigma^n$. Thus
$
	O_n(\varphi,D_\bfp)=
	\sum_{\bfw \in \Sigma^n} O(\varphi,D_{\bfp \bfw}) 
	\leq N^{2n} \lambda  \sqrt{2}|I| N^{-n-|\bfp|}
	=   \sqrt{2} \lambda  |I|N^{n-|\bfp|}.
$
\end{proof}

Now, we recall the definition of box dimension. 
For any $n_1, n_2,\ldots, n_d\in\mathbb{Z}$ and $\varepsilon>0$, we call $\Pi_{i=1}^{d}[n_i\varepsilon,(n_i+1)\varepsilon]$ an $\varepsilon$-coordinate cube in $\mathbb{R}^d$. 
Let $E$ be a bounded set in $\mathbb{R}^d$ and $\mathcal{N}_E(\varepsilon)$ the number of $\varepsilon$-coordinate cubes intersecting $E$. We define
\begin{equation}\label{eq:box-dim-def}
  \overline{\dim}_B E =\varlimsup_{\gep\rightarrow0}\frac{\log \mathcal{N}_{E}(\gep)}{\log1/\gep} \qquad \textrm{and } \qquad \underline{\dim}_B E = \varliminf_{\gep\rightarrow0}\frac{\log \mathcal{N}_{E}(\gep)}{\log1/\gep},
\end{equation}
and call them the \emph{upper box dimension} and the \emph{lower box dimension} of $E$, respectively. 
If $\overline{\dim}_B E = \underline{\dim}_B E$, we use $\dim_B E$ to denote the common value and call it the \emph{box dimension} of $E$.

Let $\varphi$ be a continuous function on $D$. 
It is well-known that in the definition of upper and lower box dimension, we can only consider $\gep_n=N^{-n}|I|$, where $n \in \mathbb{Z}^+$. That is,
\begin{equation}\label{eq:box-dim-equiv-def}
  \overline{\dim}_B \Gamma \varphi =\varlimsup_{n\to\infty} \frac{\log \N_{\Gamma \varphi}(\gep_n)}{n\log N} 
  \qquad \textrm{and } \qquad
  \underline{\dim}_B \Gamma \varphi= \varliminf_{n\to\infty}\frac{\log \N_{\Gamma \varphi}(\gep_n)}{n\log N}.
\end{equation}
It is well-known that $\underline{\dim}_B \Gamma\varphi\geq 2$ when $\Gamma \varphi$ is the graph of a continuous function on a rectangular domain of $\bR^2$. See~\cite{F90} for details.

We will use the following simple lemma in our paper, which presents a method to calculate the box dimension of the graph of a function by oscillation sums. There are results similar to this lemma in \cite{F90,JR23,RSY09}.

\begin{lemma}\label{lem:box-dim}
Let $\varphi$ be a continuous function on $D$. Then,
\[	
	1+ \varliminf_{n\to\infty}\frac{\log \big( O_n(\varphi,D)+N^n \big)}{n\log N}
	\leq \underline{\dim}_B \Gamma \varphi
	\leq \overline{\dim}_B \Gamma \varphi
	\leq 1+ \varlimsup_{n\to\infty}\frac{\log\big( O_n(\varphi,D)+N^n \big)}{n\log N}.
\]
\end{lemma}

\begin{proof}

Notice that 
\[
\N_{\Gamma \varphi}(\gep_n)	\geq  \sum_{\bfw \in \Sigma^n}  \gep_n^{-1}  O(\varphi,D_{\bfw})	\geq N^n|I|^{-1} O_n(\varphi,D).
\]
Combining this with $\N_{\Gamma \varphi}(\gep_n) \geq N^{2n}$, we have 
\begin{align*}\label{eq:box-dim-equiv-def}
	\underline{\dim}_B \Gamma \varphi
	\geq  \varliminf_{n\to\infty}   \frac{\log \big( \frac{1}{2}( N^n|I|^{-1} O_n(\varphi,D)+N^{2n} )\big)}{n\log N}   
	=1+ \varliminf_{n\to\infty}\frac{\log \big( O_n(\varphi,D)+N^n \big)}{n\log N}.
\end{align*}


On the other hand, we notice that $\N_{\Gamma \varphi}(\gep)$ and $\N_{\Gamma \varphi}(\gep_n)$ can be replaced by $\widetilde{\N}_{\Gamma \varphi}(\gep)$ and $\widetilde{\N}_{\Gamma \varphi}(\gep_n)$ in  \eqref{eq:box-dim-def} and \eqref{eq:box-dim-equiv-def} respectively, where $\widetilde{\N}_{\Gamma \varphi}(\gep)$ is the smallest number of cubes of side $\gep$ that cover $\Gamma \varphi$. Please see \cite{F90} for details. In our case,

\begin{align*}
\widetilde \N_{\Gamma \varphi}(\gep_n)\leq \sum\limits_{\bfw \in \Sigma^n} \big(\gep_n^{-1}O(\varphi,D_\bfw)+2\big) 
= N^n |I|^{-1} O_n(\varphi,D) +2N^{2n},
\end{align*}
so that
\begin{align*}
	\overline{\dim}_B \Gamma \varphi
	&=\varlimsup_{n\to\infty}\frac{\log \widetilde \N_{\Gamma \varphi}(\gep_n)}{n\log N} \\
	&\leq \varlimsup_{n\to\infty}\frac{\log \big( N^n |I|^{-1} O_n(\varphi,D) +2N^{2n}\big)   }{n\log N}
	= 1+ \varlimsup_{n\to\infty}\frac{\log \big(O_n(\varphi,D)+N^n\big)}{n\log N}.
\end{align*}
Thus, the lemma holds.
\end{proof}

\section{Vertical scaling matrices of generalized affine FISs}

In \cite{JR23,JR24}, vertical scaling matrices play the crucial role to obtain the box dimension of the graph of generalized affine FIFs in one dimensional case. 
In this section,  we will follow the spirit of these works to define vertical scaling matrices of generalized affine FISs.

\subsection{Relationship between oscillations.}

From $\Psi_w \big(x,y,f(x,y)\big)\in \Gamma f$ for all $(x,y) \in D$, we have the following useful property.
\begin{equation*}
  f\big(L_w(x,y)\big) = F_{w}\big(x,y,f(x,y)\big), \qquad \forall (x,y)\in D,\ \forall w \in \Sigma .
\end{equation*}
Thus, for any $w \in \Sigma$ and $(x,y) \in D_w$,
\begin{equation}\label{eq: FIF-rec-formula-new}
	f(x,y)=S(x,y) f\big(L_w^{-1}(x,y)\big)+q(x,y).
\end{equation}

Let $M_f:=\sup_{(x,y) \in D} |f(x,y)|$.

\begin{lemma}\label{lem:pass-1}
For any $w \in \Sigma$ and $E \subset D_w$,
\begin{align}
	 O(f,E) \leq \sup_{t \in E} |S(t)| \cdot O(f,L_w^{-1} (E))+2M_f  O(S,E)+O(q,E ),& 
	\quad \mbox{and}\label{eq:passing-1}		\\
	 O(f,E) \geq \inf_{t \in E} |S(t)| \cdot O(f,L_w^{-1} (E))-2M_f  O(S,E)-O(q,E)  .&\label{eq:passing-2}	
\end{align}
\end{lemma}
\begin{proof}
Given $ E \subset D_w$, from \eqref{eq: FIF-rec-formula-new},
\begin{align*}
	O(f,E)
	&\leq \sup_{x',x'' \in E}\big| S(x')f(L_w^{-1}(x'))-S(x'')f(L_w^{-1}(x''))\big|+\sup_{x',x'' \in E}\big|q(x')-q(x'')\big| \\
	&=\sup_{x',x'' \in E}\big| S(x')f(L_w^{-1}(x'))-S(x'')f(L_w^{-1}(x''))\big|+O(q,E).
\end{align*}
For any $x',x'',t\in E$,
\begin{align*}
	&|S(x')f(L_w^{-1}(x'))-S(x'')f(L_w^{-1}(x''))| \\
	\leq & \Big|\big(S(x')-S(t)\big) f(L_w^{-1}(x'))\Big| + \Big|\big(S(x'')-S(t)\big) f(L_w^{-1}(x''))\Big| \\
	&	   + \Big|S(t)\big(f(L_w^{-1}(x'))-f(L_w^{-1}(x''))\big)\Big| \\
	\leq & 2 M_f O(S,E)+ \sup_{t \in E} |S(t)| \cdot O(f,L_w^{-1}(E)).
\end{align*}
Thus, \eqref{eq:passing-1} holds. Using the similar argument, we can see that \eqref{eq:passing-2} holds.
\end{proof}

For any $\bfw \in \Sigma^*$, let $\overline{s}_\bfw =\sup_{x \in D_\bfw} |S(x)|$ and $\underline{s}_\bfw =\inf_{x \in D_\bfw} |S(x)|$.

\begin{lemma}\label{cor:3}
Let $\beta=\sqrt{2}(2 \lambda_S M_f + \lambda_{q} ) |I| $.
Then for any $n,k \in \mathbb{Z}^+$ and $\bfw  \in \Sigma^{n}$,
\begin{align}
	O_k(f,D_{\bfw})
	&\leq 
	\sum_{w' \in \Sigma} \overline{s}_{\bfw w'}   O_{k-1}(f,D_{\ctbfw w'}) +\beta N^{k-n} ,
	\mbox{ and} \label{eq:key-pass-11} \\
	O_k(f,D_{\bfw})
	&\geq 	\sum_{w' \in \Sigma}\underline{s}_{\bfw w'}  O_{k-1}(f,D_{\ctbfw w'}) -\beta N^{k-n}\label{eq:key-pass-12}  .
\end{align}
\end{lemma}
\begin{proof}

Fix $n,k \in \mathbb{Z}^+$ and $\bfw=w_1 \cdots w_n \in \Sigma^n$.
 For any $w' \in \Sigma$ and $\bfp \in \Sigma^{k-1}$, we have $D_{\bfw w' \bfp} \subset D_{\bfw w'} \subset D_{w_1}$ and $\ups_{\bfw w' \bfp} \leq \ups_{\bfw w'}$. Thus, from Lemma~\ref{lem:pass-1}, 
\begin{align*}	
	O(f,D_{\bfw w' \bfp}) 
	&\leq \sup_{t \in D_{\bfw w' \bfp}} |S(t)| \cdot O(f,L_{w_1}^{-1} \big( D_{\bfw w' \bfp})\big)+2M_f  O(S, D_{\bfw w' \bfp})+O(q, D_{\bfw w' \bfp})	 \\
	&\leq \ups_{\bfw w'}   O(f, D_{\sigma(\bfw )w' \bfp})+2M_f  O(S, D_{\bfw w' \bfp})+O(q, D_{\bfw w' \bfp}).	
\end{align*}
Thus,
\begin{align*}
	O_k(f,D_{\bfw}) 
	& =\sum_{w' \in \Sigma }  \sum_{\bfp \in \Sigma^{k-1}}  O(f,D_{\bfw w'\bfp})   \\
	& \leq  \sum_{w' \in \Sigma}\sum_{\bfp \in \Sigma^{k-1}} \Big( \ups_{\bfw w' } O(f,D_{\sigma(\bfw )w'\bfp}) + 2M_f  O(S,D_{\bfw w'\bfp})+O(q,D_{\bfw w'\bfp}) \Big) \\
	& = \sum_{w' \in \Sigma } \Big(  \ups_{\bfw w' } O_{k-1}(f,D_{\sigma(\bfw )w'})  \Big)+  2M_f  O_k(S,D_{\bfw })+O_k(q,D_\bfw)  .
\end{align*}
From Lemma~\ref{lem:q-lips},
\begin{equation*}
2M_f  O_k(S,D_{\bfw })+O_k(q,D_\bfw)
\leq  2 M_f \sqrt{2}\lambda_S |I| N^{k-n} + \sqrt{2}\lambda_q |I| N^{k-n}
=  \beta N^{k-n}.
\end{equation*}
Thus, \eqref{eq:key-pass-11} holds.
Similarly, we have \eqref{eq:key-pass-12} holds.
\end{proof}

\subsection{Oscillation vectors and vertical scaling matrices.}

Let $\varphi$ be a continuous function on $D$. 
Given $n,k \in \mathbb{N}$, we define a vector $V(\varphi,n,k)$ on $\Sigma^n$ by 
\[
  \big(V(\varphi,n,k)\big)_\bfj=O_k(\varphi,D_\bfj), \quad \forall \bfj \in \Sigma^n,
\]  
and call it the \emph{oscillation vector} of $\varphi$ with respect to $(n,k)$.
It is obvious that the following equality always holds,
\begin{equation}\label{eq:vector-norm}
	O_{n+k}(\varphi,D)={||{V}(\varphi,n,k)||}_1,
\end{equation}
where $||\bfx ||_1=\sum_{\bfj\in \Sigma^n} |x_\bfj|$ represents the \emph{sum norm} of $\bfx=(x_\bfj)_{\bfj\in \Sigma^n}$.

We define a $\Sigma^n \times \Sigma^n$ matrix $\upM_n$ as follows.
For any $\bfi=i_1 \cdots i_n$ and $\bfj= j_1 \cdots j_n \in \Sigma^n $, 
\begin{equation*}
	{(\upM_n)}_{\bfi,\bfj}= 
	\begin{cases}
		\overline{s}_{\bfi j_n},		&\qquad \mbox{if } \sigma(\bfi)=j_1 \cdots  j_{n-1} , \\
		0,	&\qquad \mbox{otherwise}.
	\end{cases}
\end{equation*}
We define another ${\Sigma^n} \times \Sigma^n$ matrix $\lowM_n$ by replacing $\overline{s}_{\bfi j_n}$ with $\underline{s}_{\bfi j_n}$. We call both $\upM_n$ and $\lowM_n$ \emph{vertical scaling matrices} with level $n$.

Now we are going to restate Lemma~\ref{cor:3} with matrices.
The following lemma reveals that how the vertical scaling matrices establish the relationship of oscillation vectors with different level.

\begin{lemma}\label{lem:vnp}
For all $n \in \mathbb{Z}^+$, let $\bfu_n$ be a vector on $\Sigma^n$ with $(\bfu_n)_\bfj=\sqrt{2}(2 \lambda_S M_f + \lambda_{q} ) N^{-n}|I|$ for all $\bfj\in \Sigma^n$. 
Then for any $k \in \mathbb{Z}^+$, 
\begin{equation}\label{eq:vsm-pass}
	\lowM_n V(f,n,k-1)-  N^{k} \bfu_n \leq	V(f,n,k)\leq\upM_n V(f,n,k-1)+  N^{k} \bfu_n.
\end{equation}
\end{lemma}
\begin{proof}
Let $\beta =\sqrt{2}(2 \lambda_S M_f + \lambda_{q} ) |I| $.
For any $n,k \in \mathbb{Z}^+$ and $\bfw \in \Sigma^n$, by the definition of $\upM_n$ and Lemma~\ref{cor:3},
\begin{align*}
	\big(\upM_n V(f,n,k-1)\big)_\bfw  
	=& \sum_{w' \in \Sigma} (\upM_n)_{\bfw,\sigma(\bfw)w'} V(f,n,k-1))_{\sigma(\bfw)w'} \\
	=& \sum_{w' \in \Sigma} \ups_{\bfw w'} O_{k-1}(f,D_{\sigma(\bfw)w'})  \\
	\geq&  O_k(f,D_{\bfw}) - \beta N^{k-n}
	= V(f,n,k)_\bfw - N^k (\bfu_n)_\bfw.
\end{align*}
Similarly,
$(\lowM_n V(f,n,k-1))_\bfw  
	= \sum_{w' \in \Sigma} \lows_{\bfw w'} O_{k-1}(f,D_{\sigma(\bfw)w'})  
	\leq V(f,n,k)_\bfw + N^k (\bfu_n)_\bfw$.
From the arbitrariness of $\bfw$, \eqref{eq:vsm-pass} holds for all $n,k \in \mathbb{Z}^+$.
\end{proof}

Given an $n \times n$ real matrix $A$, we define $\rho(A) = \max\{|\lambda| : \lambda \in r(A)\}$, where $r(A)$ represents the set of all eigenvalues of $A$. We call $\rho(A)$ the \emph{spectral radius} of $A$.

\begin{thm} [\cite{JR24}]\label{th:rho-rho}
For any $n \in \mathbb{Z}^+$, we have
$\rho(\lowM_n) \leq \rho(\lowM_{n+1})  \leq \rho(\upM_{n+1}) \leq \rho(\upM_n)$.
Hence, there exist $\rho_*$, $\rho^* \in \mathbb{R}$ such that $\rho^*=\lim_{n \to \infty} \rho(\upM_n)$ and $\rho_*=\lim_{n \to \infty} \rho(\lowM_n)$. Moreover, $\rho_* \leq \rho^*$.
\end{thm}

\begin{thm}[\cite{JR23}]
Assume that the vertical scaling function $S$ is either positive or negative on $D$. Then $\rho_*=\rho^*$. 
We denote the common value by $\rho_S$.
\end{thm}

A nonnegative matrix $A=(a_{ij})_{n\times n}$ is called \emph{irreducible} if for any $i,j\in \{1,\ldots,n\}$, there exists a finite sequence $i_0,\ldots,i_t\in \{1,\ldots,n\}$, such that $i_0=i,i_t=j$ and $a_{i_{\ell-1},i_\ell}>0$ for all $1\leq \ell \leq t$.
The matrix $A$ is called \emph{primitive} if there exists $k\in \bZ^+$, such that $A^k>0$. It is clear that a primitive matrix is irreducible.
Please see \cite[Chapter 8]{HJ13} for details.
The following lemma is well known. 

\begin{lemma}[Perron-Frobenius Theorem]\label{th:PF}
Let $A=(a_{ij})_{n\times n}$ be an irreducible nonnegative matrix. Then
\begin{enumerate}
	\item $\rho(A)$ is positive,
	\item $\rho(A)$ is an eigenvalue of $A$ and has a  positive eigenvector,
	\item $\rho(A)$ increases if any element of $A$ increases.
\end{enumerate}
\end{lemma}

Using the same arguements as the proof in \cite{JR23,JR24}, we have the following lemmas.
\begin{lemma}[\cite{JR23}]\label{lem:primitive-old}
If the function $S$ is not identically zero on every subrectangle of $D$, then $\upM_n$ is primitive for all $n \in \mathbb{Z}^+$.
If the function $S$ is either positive or negative on $D$, then $\lowM_n$ is primitive for all $n \in \mathbb{Z}^+$.
\end{lemma}

Define $\gamma :D \to \mathbb{R}$ by 
\begin{equation*}
	\gamma(x,y)=\sum_{w \in \Sigma}\big|S( L_w (x,y) )\big|,
\end{equation*}
which is called the \emph{sum function} with respect to $S$.
Write $\gamma_*=\min_{(x,y) \in D} \gamma(x,y) $.

\begin{lemma}\label{lem:primitive-new}
Assume that $\gamma_* \geq 1$ and the function $S$ has at most finitely many zero points.
Then there exists $n_0 \in \mathbb{Z}^+$, such that $\lowM_n$ is primitive for all $n> n_0$.

\end{lemma}

\begin{proof}
Let $Z$  be the set of the zero points of $S$.
Since $Z$ is a finite set, 
there exists $p_1 \in \mathbb{Z}^+$, such that $\sum_{w' \in \Sigma} \# ( Z \cap D_{\bfw w'} )\leq 1$ for all $\bfw \in \Sigma^n$ with $n>p_1$.  
Thus, every row of $\lowM_n$ has at least $N^2-1$ positive entries.

Let $\underline\gamma _n$ be the minimal column sum of $\lowM_n$ for all $n \in \mathbb{N}$, i.e., $\underline \gamma_n =\min_{\bfj \in \Sigma^{n}} \sum_{\bfi \in \Sigma^n} (\lowM_n)_{\bfi,\bfj} $. 
From \cite[Lemma 3.6]{JR23}, $\lim_{n \to \infty}\underline \gamma_n = \gamma_* $. 
Let $s^*:= \max\{ |S(x,y)|:(x,y) \in D \}$.
Then $s^*<1$. Thus, from $\gamma_*\geq 1$, there exists $p_2 \in \mathbb{Z}^+$, such that $\underline\gamma_n  >  s^*$ for all $n>p_2$, which implies that every column has at least $2$ positive entries.

Let $p_0=\max\{ p_1,p_2 \}$. 
Then for any $n>p_0$, every row of $\lowM_n$ has at least $N^2-1$ positive entries and every column of $\lowM_n$ has at least $2$ positive entries.
From \cite[Theorem 4.11]{JR24}, we know that there exists $n_0\in \mathbb{Z}^+$, such that $\lowM_n$ is primitive for all $n>n_0$.
\end{proof}

\section{Calculation of box dimension of fractal interpolation surfaces}

\begin{thm}\label{th:upper}
Assume that the function $S$ is not identically zero on every subrectangle of $D$. 
Then, $\overline{\dim}_{B} \Gamma f  \leq  \max\big\{2,1+\log  \rho^* /\log N \big\}$.	
\end{thm}
\begin{proof}
Fix $n \in \mathbb{Z}^+$. 
Let $\bfu_n$ be the positive constant vector mentioned in Lemma~\ref{lem:vnp}. 
From Lemma~\ref{lem:primitive-old}, $\upM_n$ is primitive.
From Lemma ~\ref{th:PF}, we can choose a strictly positive eigenvector $\bfv$ of $\upM_n$ associated with eigenvalue $\rho(\lowM_n)$  such that $\bfv\geq V(f,n,1)$ and $\bfv\geq \bfu_n$. 
Hence, from \eqref{eq:vsm-pass}, for any $p \geq 2$,
\begin{align*}
	V(f,n,p)\leq N^{p}\bfu_n+\upM_n V(f,n,p-1) \leq N^{p} \bfv+\upM_n V(f,n,p-1).
\end{align*}

From the above inequality, for all $p\geq 2$,
\begin{align*}
	V(f,n,p) \leq & N^p \bfv+\upM_n N^{p-1} \bfv+\cdots +(\upM_n)^{p-2} N^2 \bfv +(\upM_n)^{p-1} V(f,n,1) \\
	\leq & \sum_{\ell=0}^{p-1} N^{p-\ell} \rho(\upM_n)^{\ell}  \bfv 
	\leq  p N^p \max \Big\{1 ,\frac{\rho(\upM_n)^p }{N^p} \Big\}\bfv.
\end{align*}
Hence,
	$O_{n+p}(f,D)=||V(f,n,p)||_1\leq   p ||\bfv||_1 (\max\{N,\rho(\upM_n)\})^p$.
From Lemma~\ref{lem:box-dim}, 
\begin{align*}
    \overline{\dim}_B \Gamma f 
    \leq & 1+\varlimsup_{p\to\infty} \frac{\log \big(O_{n+p}(f,D)+N^{n+p}\big)}{(n+p)\log N}\\
    \leq & 1+\varlimsup_{p\to\infty} \frac{\log \big(||\bfv||_1 p (\max\{\rho (\upM_n) ,N\})^p+N^{n+p}\big)}{p\log N}
    = \max\Big\{1+ \frac{  \log \rho (\upM_n)}{\log N} ,2\Big\}.
\end{align*}

By the arbitrariness of $n$, we know from Theorem~\ref{th:rho-rho} that the theorem holds.
\end{proof}


Now, we are going to estimate the lower box dimension of generalized affine FISs satisfying one of the following two conditions.
\begin{itemize}
	\item[(A4)] $\gamma_* \geq 1$ and the function $S$ has at most finitely many zero points on $D$,
	\item[(A5)] The function $S$ is either positive or negative on $D$. 
\end{itemize}

From 	 Lemmas ~\ref{lem:primitive-old} and  ~\ref{lem:primitive-new} , if one of conditions (A4) and (A5) holds, then there exists a nonnegative integer $n_0$, such that $\lowM_n$ is primitive for all $n>n_0$. In the sequel of the paper, we always use $n_0$ to denote this constant.

\begin{lemma}\label{lem:infty}
Assume that $\varlimsup _{p \to \infty}O_p(f,D)/N^p=+\infty $ and one of the conditions (A4) and (A5) is satisfied.
Then for any $n > n_0$ and vector $\bfv$ on $\Sigma^n$, there exists $p_0 \in \mathbb{Z}^+$ such that $V(f,n,p_0) > N^{p_0}\bfv$. 
\end{lemma}

\begin{proof}



Fix $n >n_0$. 
Then $\lowM_n$ is primitive
and there exists $r \in \mathbb{Z}^+$ such that $\lowM_n^r>0$. Write $\tdM:=N^{-r}\lowM_n^r$ and $\delta=\min_{\bfi,\bfj \in \Sigma^n}\tdM_{\bfi,\bfj}$ be the smallest entry of $\tdM$. It is clear that $\delta>0$.

Let $\bfu_n$ be vector mentioned in Lemma~\ref{lem:vnp}.
From \eqref{eq:vsm-pass}, for any $p \in \mathbb{Z}^+$, we have
\begin{align*}
	V(f,n,p+r)\geq \lowM_n^r V(f,n,p)- N^{p+r}\bfu_n - N^{p+r-1} \lowM_n \bfu_n  -\cdots- N^{p+1} \lowM_n^{r-1} \bfu_n .
\end{align*}
Write $\tdu:=\sum_{\ell=0}^{r-1}  N^{-\ell} \lowM_n^{\ell} \bfu_n   $.  
Notice that the above inequality can be rewritten as
\begin{align}\label{eq:v-pass-org}
	V(f,n,p+r)\geq N^r \tdM V(f,n,p)- N^{p+r} \tdu.
\end{align}
Divide \eqref{eq:v-pass-org} by $N^{p+r}$, for any $\bfw \in \Sigma^n$,
\begin{align*}
		\frac{V(f,n,p+r)_\bfw}{N^{p+r}} 
	\geq 
	 \sum_{\bfi \in \Sigma^n}\tdM_{\bfw,\bfi} \frac{ V(f,n,p)_\bfi}{N^p}- \tdu_\bfw
	\geq \delta \frac{ O_{n+p}(f,D)}{N^p}-  \tdu_\bfw.
\end{align*}

Fix $\bfv$ on $\Sigma^n$.
Let $G=\max\{ (\tdu_\bfk+\bfv_\bfk) :\bfk \in \Sigma^n\} $. 
From  $\varlimsup _{p \to \infty}O_p(f,D)/N^p=+\infty $, there exists $p_1 \in \mathbb{Z}^+$ satisfying $O_{p_1}(f,D)/ N^{p_1}>G / \delta$, which implies
\begin{align*}
		\frac{V(f,n,p_1+r)_\bfw}{N^{p_1+r}} 
	\geq \delta \frac{ O_{n+p_1}(f,D)}{N^{p_1}}-  \tdu_\bfw
	\geq \delta \frac{ O_{p_1}(f,D)}{N^{p_1}}-  \tdu_\bfw
	>G -\tdu_\bfw \geq 	\bfv_\bfw.
\end{align*}
Let $p_0=p_1+r$. 
Then $V(f,n,p_0)>N^{p_0}\bfv$.
\end{proof}

\begin{thm}\label{th:lower}
Assume that one of the conditions (A4) and (A5) is satisfied.
Then in the case that $\varlimsup _{p \to \infty}O_p(f,D)/N^p=+\infty $, we have 
$\underline{\dim}_{B} \Gamma f  \geq  1+ \log  \rho_* /\log N$.
\end{thm}
\begin{proof}

If $\rho_* \leq N$, then the lemma always holds because $\underline{\dim}_B \Gamma f \geq 2$. We assume that $\rho_*>N$ in the following part.

Fix $n >n_0$ satisfying $\rho(\lowM_n)>N$. Then $\lowM_n$ is primitive. 
Let $\bfu_n$ be the constant positive vector mentioned in Lemma~\ref{lem:vnp}.
From \eqref{eq:vsm-pass}, for all $p \in \mathbb{Z}^+$,
\begin{equation}\label{eq:tau}
	\frac{V(f,n,p)}{N^p}\geq \frac{ \lowM_n }{N} \frac{V(f,n,p-1)}{N^{p-1}} -\bfu_n.
\end{equation}

For any $N<\tau<\rho(\lowM_n) $, from Lemma ~\ref{th:PF}, we can choose a positive eigenvector $\bfv$ of $\lowM_n$ associated with eigenvalue $\rho(\lowM_n)$ satifying $\bfv \geq N \bfu_n/(\rho(\lowM_n)-\tau)$. Combining this with Lemma~\ref{lem:infty}, there exists $ p_0 \in \mathbb{Z}^+$ such that 
\begin{align*}
	\frac{V(f,n,p_0)}{N^{p_0}}\geq \bfv \geq \frac{N}{\rho(\lowM_n)-\tau}\bfu_n.
\end{align*}
Notice that 
\[
  \frac{ \lowM_n }{N}\bfv  - \bfu_n \geq \frac{\rho(\lowM_n)}{N} \bfv - \frac{\rho(\lowM_n)-\tau}{N} \bfv =\frac{\tau}{N} \bfv.
\]
Hence, from \eqref{eq:tau}, 
\begin{equation*}
	\frac{V(f,n,p_0+1)}{N^{p_0+1}}\geq \frac{ \lowM_n }{N} \frac{V(f,n,p_0)}{N^{p_0}} -\bfu_n
	\geq \frac{ \lowM_n }{N}\bfv  -\bfu_n
	=\frac{\tau}{N} \bfv.
\end{equation*}
By induction, we can prove that for all $p \in \mathbb{Z}^+$, 
\begin{equation*}
	\frac{V(f,n,p_0+p)}{N^{p_0+p}}
	\geq \frac{ \lowM_n }{N} \Big(\frac{\tau^{p-1}}{N^{p-1}}-1 \Big) \bfv 
	 +\frac{\lowM_n }{N} \bfv-\bfu_n
	\geq \frac{\tau}{N} \Big(\frac{\tau^{p-1}}{N^{p-1}}-1\Big) \bfv +\frac{\tau}{N} \bfv = \frac{\tau^p}{N^p} \bfv.
\end{equation*} 
Hence, 
$V(f,n,p_0+p) \geq \tau^pN^{p_0} \bfv.$ 
From \eqref{eq:vector-norm}, we have $O_{n+p_0+p}(f,D)\geq \tau^p N^{p_0} ||\bfv||_1 $. 
Thus, from Lemma ~\ref{lem:box-dim}, 
\begin{align*}
\underline{\dim}_{B} \Gamma f  \geq 1+\varliminf_{k \to \infty} \frac{\log (O_k(f,D)+N^k)}{k \log N}	
	  \geq 1+\varliminf_{p \to \infty} \frac{\log (\tau^p N^{p_0} ||\bfv||_1)}{(n+p_0+p) \log N} 
	  = 1+\frac{\log \tau}{ \log N}.
\end{align*}

From the arbitrariness of $\tau$,
$
	\underline{\dim}_{B} \Gamma f  
	\geq  1+\log \rho(\lowM_n) / \log N.
$
Combining this with Theorem~\ref{th:rho-rho}, we have
	$\underline{\dim}_{B} \Gamma f  
	\geq 1+\log \rho_* / \log N.$
\end{proof}

\begin{thm}\label{th:positive}
Assume that the function $S$ is either positive or negative on $D$. 
Then in the case that $\rho_S>N$ and $\varlimsup _{p \to \infty}O_p(f,D)/N^p=+\infty $,
\begin{equation*}
	{\dim}_{B} \Gamma f  =  1+\frac{\log  \rho_S   }{\log N}.
\end{equation*}
Otherwise, $ {\dim}_{B} \Gamma f  =2$.
\end{thm}
\begin{proof}

Since the function $S$ is either positive or negative, from Lemma ~\ref{lem:primitive-old}, both $\upM_n$ and $\lowM_n$ are primitive for all $n \in  \mathbb{Z}^+$.

If there exists $c \in \mathbb{R}$ satisfying $\varlimsup _{p \to \infty}O_p(f,D)/N^p=c$, then from Lemma ~\ref{lem:box-dim}, 
\[
	\overline\dim_B \Gamma f \leq 1+ \varlimsup_{ n \to \infty} \frac{\log (O_n(f,D)+N^n)}{n \log N}=2,
\]
so that $ {\dim}_{B} \Gamma f  =2$.

In the case that $\rho_S \leq N$, from Theorem~\ref{th:upper}, 
$
	\overline\dim_B \Gamma f \leq \max \big\{ 2,1+ \log \rho_S / \log N \big\} =2.
$
Thus, $ {\dim}_{B} \Gamma f  =2$.

In the case that $\rho_S>N$ and $\varlimsup _{p \to \infty}O_p(f,D)/N^p=+\infty $, we know from Theorems \ref{th:upper} and \ref{th:lower} 
that ${\dim}_{B} \Gamma f  =  1+ \log  \rho_S / \log N.$ 
\end{proof}

It is not intuitive to check the condition $\varlimsup_{p \to \infty} O_p(f,D)/N^p =+\infty$ mentioned in Theorem \ref{th:positive} directly. 
Now we will provide a criterion.

Let $\beta=\sqrt{2}(2 \lambda_S M_f + \lambda_{q} ) |I| $.
Using the similar arguement as Lemma ~\ref{lem:pass-1}, for any $k \in \mathbb{Z}^+$, we have
\begin{equation}\label{eq:gamma-pass}
	O_k (f,D) 	 \geq  \gamma_*  O_{k-1}(f,D)-   \beta N^k.
\end{equation}
Assume that $\gamma_* \neq N$.
Let $c={\beta N}/{(\gamma_*-N)}$.
Then \eqref{eq:gamma-pass} can be rewritten as
\[
	\frac{O_k(f,D)}{N^k}-c \geq \frac{\gamma_*}{N}\Big(\frac{O_{k-1}(f,D)}{N^{k-1}}-c \Big).
\]
Thus, if $\gamma_*>N$ and $O_{k} (f,D)/N^{k}>c$ holds for some $k \in \mathbb{Z}^+$, 
then $\varlimsup_{p \to \infty} O_p(f,D)/N^p =+\infty$.

\ack 
The author wishes to thank Professor Huo-Jun Ruan for many helpful discussions and valuable suggestions. 
The author also wishes to thank Doctor Zhen Liang and Doctor Jian-Ci Xiao for their helpful advices.


\end{document}